\documentclass[12pt]{article}
\usepackage{amsmath}
\usepackage{amssymb}
\usepackage{verbatim}

\newcommand{\syl}[1]{(\mkern -5mu(\mathbf{#1})\mkern -5mu)}
\newcommand{\den}[1]{\frac{1}{1-x^{#1}}}
\newcommand{\dena}[1]{\frac{1}{1-ax^{#1}}}
\newcommand{\ven}[1]{1-x^{#1}}

\newcommand{\nena}[1]{1+ax^{#1}}
\newcommand{\gauss}[2]{\mathbf{\binom{#1}{#2}}}

\begin{document}
\title{ TWO FORMULAS FOR  PLANE PARTITIONS}
\author{Peter Lindqvist}
\date{\footnotesize{Department of Mathematical Sciences\\Norwegian University of Science and Technology\\NO--7491 Trondheim, Norway}}
\maketitle

\begin{center}
  \texttt{Dedicated to the memory of Jaak Peetre 1935-2019}
\end{center}
 \bigskip
{\small \textsc{Abstract:}}\footnote{AMS classification 05A17, 11P81} \textsf{We give a direkt deduction and proof of two identities in the theory of partitions. The first one is known to enumerate the traces of plane partitions. The second one is not  given any combinatorial interpretation.}

\section{Introduction}
The generating formula 
\begin{equation}
  \den{1}\,\den{2}\,\den{3}\,\cdots\,= \,1 + \sum_{n=1}^{\infty}p(n)x^n
\end{equation}
for partitions was given by Euler.
For example, there are $p(5) = 7$ partitions of the integer $5$:
\begin{gather*}
  5=4\!+\!1=3\!+\!2=3\!+\!1\!+\!1=2\!+\!2\!+\!1=2\!+\!1\!+\!1\!+\!1=1\!+\!1\!+\!1\!+\!1\!+\!1.
\end{gather*}
 The identity
\begin{gather}\label{Euler}
  \dena{1}\,\dena{2}\,\dena{3}\,\cdots\,\\=\, 1+\den{1}xa+\den{1}\,\den{2}x^2a^2+\den{1}\,\den{2}\,\den{3}x^3a^3+\cdots \nonumber
\end{gather}
appeared 1748 in Euler's \emph{Introductio in Analysin Infinitorum}, see  paragraph  313, Chapter XVI of \cite{E}. This \footnote{Equation (2.2.5) in \cite{A} suffers from a misprint. This general introduction to partitions  also includes plane partitions.} is  Theorems 349 and 350 in  \cite{H} or formula [5b] in section 2.5 of \cite{C}.

P. MacMahon \cite{M} derived the generating formula
\begin{equation}
  \frac{1}{1-x}\,\frac{1}{(1-x^2)^2}\,\frac{1}{(1-x)^3}\,\cdots\, =\, 1 + \sum_{n=1}^{\infty}p\mkern - 2mu p(n)x^n
\end{equation}
for \emph{plane partitions}. The expansion begins
$$1+x+3x^2+6x^3+13x^4+24x^5+\cdots.$$
{\small In passing, we cite the asymptotic formula (for large $n$)
\begin{equation*}
  p\mkern -2mup(n)\,\approx \, \frac{\zeta(3)^{\frac{7}{36}}}{\sqrt{12\pi}}\Bigl(\frac{n}{2}\Bigr)^{\frac{25}{36}}\exp \Bigl(3\,\zeta(3)^{\frac{1}{3}}\Bigl(\frac{n}{2}\Bigr)^{\frac{2}{3}} + \zeta'(-1)\Bigr)
\end{equation*}
due to  E. Wright, cf. \cite{W}.}

The counterpart to Euler's formula (\ref{Euler}) is 
\begin{align}\label{Stanley}
  &\frac{1}{1-ax}\,\frac{1}{(1-ax^2)^2}\,\frac{1}{(1-ax)^3}\,\cdots\\&\, =\,1 + \sum_{n=1}^{\infty}\frac{g_n(x)}{(1-x)^2(1-x^2)^2(1-x^3)^2 \cdots(1-x^n)^2}\,\,x^na^n,\nonumber
\end{align}
where the $g_n(x)$ are polynomials with positive integer coefficients.  R. Stanley found this formula and gave a combinatorial interpretation for the coefficients, cf. \cite{S}. The expansion begins
\begin{gather*}
1+\frac{1}{[\ven{1}]^2}\,xa+\frac{1+x^2}{[(\ven{1})(\ven{2})]^2}\,x^2a^2\\+\frac{1+x^2+2x^3+x^4+x^6}{[(\ven{1})(\ven{2})(\ven{3})]^2}\,x^3a^3+\cdots \nonumber
\end{gather*}

The object of our note is first to provide a very simple derivation of Stanley's formula, using the same method as Euler used for his formula in  \cite{E}. There is hardly anything new here. The proof does not shed light on the combinatorial interpretation with traces. We avoid the exponentiation with  Young polynomials and Bruno di Fa\`{a}'s formula, but the Gaussian polynomials
$$\binom{\mathbf{n}}{\mathbf{k}} =\frac{(1-x)^{n-k+1}(1-x)^{n-k+2}\cdots(1-x^n) }{(1-x)(1-x^2)\cdots(1-x^k)},\quad k=0,1,2,...,n,$$
are intrinsic also in our  proof. \emph{They are polynomials in $x$ with positive integer coefficients}. As $x \to 1$, they become the usual binomial coefficients. See paragraph 110 in \cite{R} for the Gaussian polynomials.

Second, we shall derive the formula
\begin{align}\label{new}
  (\nena{1})&(\nena{2})^2(\nena{3})^3\,\cdots\,\\
  =\, 1& + \sum_{n=1}^{\infty}\frac{h_n(x)}{(\ven{1})^2(\ven{2})^2(\ven{3})^3 \cdots(\ven{n})^2}\,x^na^n, \nonumber
\end{align}
where the $h_n(x)$ are polynomials.

\paragraph{Notation.}
The notation of Cayley is convenient. He wrote
$$\syl{m} = 1-x^m.$$
In this notation
$$\syl{1}\syl{2}\syl{3}\cdots\syl{m} = (1-x)(1-x^2)(1-x^3)\cdots (1-x^m).$$ In particular, (\ref{Stanley}) becomes
$$\boxed{\prod_{n=1}^{\infty}(1-ax^n)^{-n}\,=\, 1 + \sum_{n=1}^{\infty}\frac{g_n(x)\,x^n}{\{\syl{1}\syl{2}\syl{3}\cdots\syl{n}\}^2}\,a^n.}$$
We can write (\ref{new}) as
$$\boxed{\prod_{n=1}^{\infty}(\nena{n})^n\,=\,1 + \sum_{n=1}^{\infty}\frac{h_n(x)\,x^n}{\{\syl{1}\syl{2}\syl{3}\cdots\syl{n}\}^2}\,a^n.}$$
\section{Preliminaries}

\emph{A plane partition} of an integer $n\geq 1$ is an array of positive intgers $n_{ij}$ with non-increasing rows and columns
 
\begin{eqnarray*}
n_{11}\quad n_{12}\quad n_{13}\quad\cdots\\
n_{21}\quad n_{22}\quad n_{23}\quad\cdots\\
n_{31}\quad n_{32}\quad n_{33}\quad\cdots\\
\,\,\vdots\qquad\vdots \qquad\vdots \qquad\ddots   \\
\end{eqnarray*}
that sum up to $n = \sum_{i,j}n_{ij}$. It contains only finitely many entries. The \emph{trace} of the partition is the diagonal sum
$$n_{11}+n_{22}+n_{33}+\cdots,$$
where we stop the summation just before the first undefined diagonal element.

Let us take a simple example. Take $n = 4$. There are $p \mkern -2mu p(4) = 13$ plane partitions of the integer $4$. They are
\begin{align*}
  4\qquad3\,1\qquad\mathbf{2}\,2\qquad\mathbf{2}\,1\,1\qquad 1\,1\,1\,1\phantom{0}\,
\end{align*}
\begin{eqnarray*}
  3\qquad\mathbf{2}\,\phantom{0}\qquad \mathbf{2}\,1 \qquad \mathbf{1}\,1\,\phantom{0} \qquad 1\,1\,1\\
  1\qquad2\,\phantom{0}\qquad 1\,\phantom{0} \qquad 1\,\mathbf{1} \,\phantom{0}\qquad 1\,\phantom{0}\,\phantom{0}\\
\end{eqnarray*}
\vspace*{-3\belowdisplayskip}
\begin{eqnarray*}
  1 \, 1\qquad\mathbf{2}\qquad 1\\
  1\, \phantom{0} \qquad 1\qquad 1\\
  1\,\phantom{1} \qquad 1\qquad 1\\
  \phantom{0}\,\phantom{0} \qquad \phantom{0} \qquad 1
\end{eqnarray*}
Among these there are $6$ having trace = $2$. The interpretation of formula (\ref{Stanley}) is that $6$ is the coefficient $c_{42}$ of $x^4a^2$ in the expansion
$$\prod_{n=1}^{\infty}\bigl(1-ax^n\bigr)^{-n}\,=\, 1 + \sum_{i,j=1}^{\infty}c_{ij}x^ia^j.$$
In particular, in (\ref{Stanley}) the $a^2$-term
$$\frac{(1+x^2)x^2}{[(1-x)(1-x^2)]^2}\,a^2\,=\,\left( x^2 +2x^3+\mathbf{6}x^4+10x^5+\cdots\right)a^2$$
 is the generating function for those plane partitions that have trace $=\,2$.

\section{The Proof for $\prod(1-ax^n)^{-n}$}

Write
$$\frac{1}{1-ax}\,\frac{1}{(1-ax^2)^2}\,\frac{1}{(1-ax^3)^3}\,\cdots\,=\, 1+c_1a+c_2a^2+c_3a^3+\cdots,$$
where $c_n=c_n(x)$ are the yet undetermined coefficients. Replace $a$ by $ax$ to obtain
$$\frac{1}{1-ax^2}\,\frac{1}{(1-ax^3)^2}\,\frac{1}{(1-ax^4)^3}\,\cdots\,=\, 1+c_1ax+c_2(ax)^2+c_3(ax)^3+\cdots.$$
Divide this formula by $(1-ax)(1-ax^2)(1-ax^3)\cdots$ to infer that
$$ 1+c_1a+c_2a^2+c_3a^3+\cdots=\frac{1+c_1ax+c_2(ax)^2+c_3(ax)^3+\cdots}{(1-ax)(1-ax^2)(1-ax^3)\cdots}.$$
Now Euler's identity (\ref{Euler}) yields
\begin{align*} 1+c_1a+c_2a^2+c_3a^3+\cdots\qquad \qquad \qquad \qquad\qquad\qquad \qquad \qquad\qquad \\
  =
  \Bigl\{1+\frac{xa}{1-x}+\frac{x^2a^2}{(1-x)(1-x^2)}+\frac{x^3a^3}{(1-x)(1-x^2)(1-x^3)}+\cdots \Bigr\}&\\
  \times\Bigl\{1+c_1ax+c_2(ax)^2+c_3(ax)^3+\cdots\Bigr\}&.
\end{align*}

The next step is to multiply the two series upon which a comparison of the coefficients yields
\begin{align*}
  c_na^n\,=\,c_na^nx^n &+ c_{n-1}a^{n-1}x^{n-1}\den{1}xa
  +c_{n-2}a^{n-2}x^{n-2}\den{1}\den{2}x^2a^2\\
  &+ c_{n-3}a^{n-3}x^{n-3}\den{1}\den{2}\den{3}x^3a^3+\cdots\cdot\cdot\\
  &+c_{n-k}a^{n-k}x^{n-k}\den{1}\den{2}\den{3}\cdots\den{k}x^ka^k +\cdots\cdot\cdot\\&+ c_1 ax\den{1}\den{2}\den{3}\cdots\den{n-1}x^{n-1}a^{n-1} \\&+ \den{1}\den{2}\den{3}\cdots\den{n}x^{n}a^{n}.
\end{align*}
Each term to the right has the factor $x^n$. We see that in (\ref{Stanley})
$$g_n(x)x^n\,=\,\bigl[(1-x)(1-x^2)(1-x^3)\cdots(1-x^n)\bigr]^2c_n\,\equiv\,b_n.$$
For $b_n = g_n(x)x^n$ it follows that

\begin{align*}
  b_nx^{-n}\,=\, b_{n-1}\frac{\syl{n}}{\syl{1}} &+ b_{n-2}\frac{\syl{n}\cdot\syl{n-1}^2}{\syl{1}\syl{2}} 
+ b_{n-3}\frac{\syl{n}\cdot\syl{n-1}^2\syl{n-2}^2}{\syl{1}\syl{2}\syl{3}}\\&+\cdots\cdot\cdot+\\ &+b_{n-k}\frac{\syl{n}\cdot\syl{n-1}^2\syl{n-2}^2\cdots \syl{n-k+1}^2}{\syl{1}\syl{2}\syl{3}+\cdots +\syl{k}}\\&+\cdots\cdot\cdot+\\ &+b_1\frac{\syl{n}\cdot\syl{n-1}^2\syl{n-2}^2\cdots\syl{2}^2}{\syl{1}\syl{2}\syl{3}+\cdots +\syl{n-1}}\\
    &+    \frac{\syl{n}\cdot\syl{n-1}^2\syl{n-2}^2\cdots\syl{2}^2\syl{1}^2}{\syl{1}\syl{2}\syl{3}+\cdots +\syl{n}}.
\end{align*}

This can succinctly be written as
 \begin{align*}
   b_nx^{-n}\,&=\, b_{n-1}\binom{\mathbf{n}}{\mathbf{1}} + b_{n-2}\binom{\mathbf{n}}{\mathbf{2}}\cdot\syl{n-1}
+ b_{n-3}\binom{\mathbf{n}}{\mathbf{3}}\cdot\syl{n-1}\syl{n-2}\\&+\cdots\cdot\cdot +b_{n-k}\binom{\mathbf{n}}{\mathbf{k}}\cdot\syl{n-1}\syl{n-2}\cdots \syl{n-k+1}\\&+\cdots\cdot\cdot
    +b_1\binom{\mathbf{n}}{\mathbf{n-1}}\cdot\syl{n-1}\syl{n-2}\cdots\syl{2}\\
    &+    \binom{\mathbf{n}}{\mathbf{n}}\cdot\syl{n-1}\syl{n-2}\cdots\syl{2}\syl{1}.
\end{align*}
 In terms of the original $g_n$ the recursion reads
 $$\boxed{g_n(x)\,=\, \sum_{k=1}^{n} \binom{\mathbf{n}}{\mathbf{k}}\syl{n-k+1}\cdots \syl{n-2}\syl{n-1}\,g_{n-k}(x) x^{n-k} .}$$
 Here $g_0(x) =1$.
 Since the Gaussian polynomials  have positive integer coefficients, so does $g_n(x)$ by induction. The  desired result follows.  \qquad $\Box$

 \section{Expansion of $\prod \bigl(1+ax^n\bigr)^n$}

 Write
 $$(\nena{1})(\nena{2})^2(\nena{3})^3\,\cdots=\, 1+d_1a+d_2a^2+d_3a^3+\cdots$$
 where $d_n=d_n(x)$. Replace $a$ by $ax$ to obtain
 $$(\nena{2})(\nena{3})^2 (\nena{4})^3 \,\cdots=\,1+d_1ax+d_2(ax)^2+d_3(ax)^3+\cdots.$$
 Upon multiplication with $(\nena{1})(\nena{2})(\nena{3})\,\cdots$ we arrive at
 \begin{equation*}
   1+d_1a+d_2a^2+\cdots = \bigl\{(\nena{1})(\nena{2})(\nena{3})\cdots\bigr\}\bigl\{1+d_1ax+d_2a^2x^2 +\cdots\bigr\}.
 \end{equation*}
 Using Euler's formula
 $$\prod_{n=1}^{\infty} \bigl(1+ax^n\bigr) \,=\, 1 + \sum_{n=1}^{\infty}\frac{x}{\ven{1}}\frac{x^2}{\ven{2}}\cdots \frac{x^n}{\ven{n}}\,a^n$$
 (Paragraph 307 in Chapter XVI of \cite{E}; Theorem 349 in \cite{H}; formula [5k] in Section 2.5 of \cite{C}) and multipying we get the recursion
 \begin{align*}
   d_n\,=\,&\frac{x}{\ven{1}}\frac{x^2}{\ven{2}}\cdots \frac{x^n}{\ven{n}}+\frac{x}{\ven{1}}\frac{x^2}{\ven{2}}\cdots \frac{x^{n-1}}{\ven{n-1}}\,d_1x\\
   +&\frac{x}{\ven{1}}\frac{x^2}{\ven{2}}\cdots \frac{x^{n-2}}{\ven{n-2}}\,d_2x^2+\,\cdots\,+\frac{x}{\ven{1}}d_{n-1}x^{n-1} + d_nx^n.
 \end{align*}
 Move the term  $d_nx^n$ to the left side, multiply by $(1-x)(1-x^2)\cdots(1-x^{n})$, and abbreviate
 $$e_n = (1-x)(1-x^2)\cdots(1-x^{n})d_n$$
 in order to  obtain
 \begin{align*}
   (1-x^n)e_n\,= &\,x^{n-1}d_{n-1}\,x\gauss{n}{1} + \,x^{n-2}d_{n-2}\,x^{1+2}\gauss{n}{2} \\
   &+ \,x^{n-3}d_{n-3}\,x^{1+2+3}\gauss{n}{3}\,+\cdots+\, 1\,x^{\frac{n(n+1)}{2}}\gauss{n}{n}
 \end{align*}
 Multiplying by $(1-x)(1-x^2)\cdots(1-x^{n-1})$ and noticing that\\ $h_nx^n = (1-x)(1-x^2)\cdots(1-x^{n})e_n$ we arrive at the final recursion
 $$\boxed{h_n(x)\,=\, \sum_{k=1}^{n} \gauss{n}{k}\syl{n-1}\syl{n-2}\cdots \syl{n-k+1}\,h_{n-k}(x)x^{\frac{k(k-1)}{2}},}$$
 where we have put $h_0=1$.\qquad $\Box$

 \bigskip

 {\small \textsf{Remark}: The expansion begins
   \begin{gather*} 1+\frac{1}{(\ven{1})^2}ax+ \frac{2x}{(\ven{1})^2(\ven{2})^2}x^2a^2+
     \frac{x^2(1+4x+x^2)}{(\ven{1})^2(\ven{2})^2(\ven{3})^2}x^3a^3\\ +
     \frac{x^4(3+4x+10x^2+4x^3+3x^4)}{(\ven{1})^2(\ven{2})^2(\ven{3})^2(\ven{4})^2  }x^4a^4\\+\frac{x^6(3+8x+15x^2+20x^3+28x^4+20x^5+15x^6+8x^7+3x^8)}{(\ven{1})^2(\ven{2})^2(\ven{3})^2(\ven{4})^2(\ven{5})^2  }x^5a^5  + \,\cdots
     \end{gather*}}


\begin{thebibliography}{A}{\small
  \bibitem[A]{A} \textsc{G. Andrews}, {\it The Theory of Partitions}, Cambridge University Press, Cambridge 1984.
    \bibitem[C]{C} \textsc{L. Comtet}, {\it Advanced Combinatorics}, D. Reidel Publishing Company, Dordrecht 1974. 
     \bibitem[HW]{H} \textsc{G. Hardy,  E. Wright}, {\it An Introduction to the Theory of Numbers}, $5^{th}$ Edition, Clarendon Press, Oxford 1979.
    \bibitem[E]{E} \textsc{L. Euler}, {\it Introductio in Analysin Infinitorum}, Lausanne 1748. English translation: {\it Introduction to Analysis of the Infinite}, Springer, New York 1988.
    \bibitem[M]{M} \textsc{P. MacMahon}, {\it Combinatory Analysis}. Volume II, Cambridge 1916.
    \bibitem[R]{R} \textsc{H. Rademacher}, {\it Topics in Analytic Number Theory}, Die Grundlehren der mathematischen Wissenschaften in Einzeldarstellungen \textbf{169}, Springer-Verlag, Berlin 1973.  
    \bibitem[S]{S} \textsc{R. Stanley}, {\it The conjugate trace and trace of a plane partition}, Journal of Combinatorial Theory (A) \textbf{14}, 1973, pp. 53--65.
    \bibitem[W]{W} \textsc{E. Wright}, {\it   Asymptotic partition formulae I. Plane partitions}, The Quarterly Journal of Mathematics (Oxford Series) \textbf{2}, 1931.
      }
\end{thebibliography}
\end{document}